\documentclass[10pt]{article}
\usepackage{amssymb,amsthm}
\usepackage[tbtags]{amsmath}
\usepackage{color}

\definecolor{c20}{rgb}{0.,0.7,0.}
\definecolor{c30}{rgb}{0.,0.,1.}
\definecolor{c40}{rgb}{1,0.1,0.7}
\definecolor{c50}{rgb}{1,0,0}
\def\IF{\infty}

\def\bE#1{\textcolor{c20}{#1}}
\def\bE#1{#1}
\def\rE#1{\textcolor{c20}{#1}}
\def\rE#1{#1}

\def\cE#1{\textcolor{c20}{#1}}
\def\cE#1{#1}

\def\cEE#1{\textcolor{c20}{#1}}
\def\cEE#1{#1}

\def\wwF{\cE{\widehat f}}

\newcommand{\COM}[1]{}
\def\QED{$\hfill \Box$\\}

\newtheorem{thm}{Theorem}[section]
\newtheorem{cor}{Corollary}[section]
\newtheorem{lem}{Lemma}[section]

\theoremstyle{definition}
\newtheorem{defn}{Definition}[section]
\newtheorem{exm}{Example}[section]
\theoremstyle{remark}

\newcommand{\BQN}{\begin{eqnarray}}
\newcommand{\EQN}{\end{eqnarray}}
\newcommand{\BQNY}{\begin{eqnarray*}}
\newcommand{\EQNY}{\end{eqnarray*}}

\newcommand{\R}{\mathbb{R}}
\newcommand{\1}{\mathbb{I}}
\newcommand{\I}{\mathbf{I}}

\newcommand{\F}{\mathcal{F}}
\def\Lim{\lim\limits}
\def\Inf{\inf\limits}
\newcommand{\E}[1]{\mathbb{E}\left\{#1 \right\}}

\newcommand{\pk}[1]{P(#1) }
\newcommand{\norm}[1]{\lVert #1 \rVert}
\def\kH{\rE{\mathcal{H}}}
\def\wF{{\wwF}}

\def\Nhh{\rE{{\mathcal{H}}}}
\def\HhH{\rE{{\mathcal{H}}}}

\begin{document}

\centerline{\bf \large Boundary non-crossing probabilities for fractional Brownian motion  with trend}

\centerline{Enkelejd  Hashorva\footnote{Department of Actuarial Science, University of Lausanne, UNIL-Dorigny, 1015 Lausanne, Switzerland,
email:enkelejd.hashorva@unil.ch},  
Yuliya Mishura\footnote{Department of Probability, Statistics and Actuarial Mathematics, National Taras Shevchenko University of Kyiv, 01601 Volodymyrska 64, Kyiv, Ukraine, email: myus@univ.kiev.ua },
and Oleg Seleznjev\footnote{             
Department of Mathematics and Mathematical Statistics,
Ume{\aa} University, SE-901 87 Ume\aa, Sweden, email:oleg.seleznjev@matstat.umu.se}}

\bigskip

{\bf Abstract}: In  this paper we investigate the boundary non-crossing probabilities of a fractional Brownian motion considering some general deterministic trend function. We derive bounds for non-crossing probabilities and discuss the case of a large trend function. As a by-product we solve a minimization problem related to the norm of the trend function.

{\bf Key Words}:  boundary crossings; Cameron-Martin-Girsanov theorem; reproducing kernel Hilbert space; large deviation principle; Molchan martingale; fractional Brownian motion.



\section{Introduction}
Calculation of boundary crossing (or non-crossing) probabilities of Gaussian processes with trend is a long-established   
 and interesting topic of applied probability, see, e.g., 
\cite{Durb92,Wang97,Nov99,MR1787122, Wang2001,MR2009980,BNovikov2005,MR2028621,1103.60040,1137.60023,1079.62047,MR2443083,MR2576883,Janssen08} and references therein. Numerous applications concerned with the evaluation of boundary non-crossing probabilities relate  to mathematical finance, risk theory, queueing theory, statistics, physics, biology among many other fields.
In the literature,  
 most of contributions treat the case when the Gaussian process $X(t),t \ge 0$ is a Brownian motion which allows to calculate the boundary non-crossing probability  $\pk{X(t)+ f(t) < u, t\in [0,T]}$,  
 for some trend function $f$ and  two  given constants $T,u>0$  
  by various methods (see, e.g., \cite{Alil2010,MR2752890}).  
   For particular $f$ including the case of a piecewise  
    constant function,  
     explicit calculations are possible, see, e.g., \cite{MR2065480}.   
     Those explicit calculations allow then to approximate the non-crossing probabilities for general $f$ and for $f$ being large, see \cite{MR2065480,MR2175400,MR2028621}.\\
In this paper the centered Gaussian process $X=B^H$ is a fractional Brownian motion (fBm) with Hurst index $H\in (0,1)$ for which no explicit calculations of the boundary non-crossing probability
are possible for the most of the trend functions.\\
 Therefore, our interest in this paper is on the derivation of upper and lower bounds for
 $$P_f:=\pk{B^H(t) + f(t)\le u(t), t \in \R_+}$$     
  for some admissible trend functions $f$ and measurable functions $u: \R_+ \to \R$ such that $u(0)\ge0$.
  In the following we shall consider \cEE{$f\not=0$} to belong to the reproducing kernel Hilbert Space (RKHS) of $B^H$ which is denoted by $\kH$ defined by the covariance kernel of $B^H$ given as
\BQN \label{Rh}
\bE{R_{H}(s,t)}:=\E{B^H(s)B^H(t)}=\frac{1}{2}(t^{2H}+s^{2H}-|t-s|^{2H}),\quad t,s\geq0.
\EQN
A precise description of  $ \HhH $ is given in Section 2,  
 where also the norm $\norm{f}_{ \Nhh }$ for $f\in  \HhH $ is defined;   \rE{for notational simplicity we suppress the Hurst index $H$ and the specification of $\R_+$ avoiding the more common notation $\mathcal{H}_H(\R_+)$}.

The lack of explicit formulas (apart from $H=1/2$ case) for trend functions $f$ and given $u$ poses problems for judging the
accuracy of our bounds for $P_f$. A remedy for that is to consider the asymptotic performance of the bounds for trend functions $\gamma f$
with $\gamma \to 0$ and $\gamma \to \IF$.
The latter case is more tractable since \cEE{if for some $x_0$ we have $f(x_0)>0$, then (see Corollary \ref{corLDP} below)}
\BQN\label{LD}
\cEE{\ln P_{\gamma f} \ge -(1+o(1))\frac{\gamma^2}{2}\norm{\wF}^2_{ \Nhh }, \quad \gamma \to \IF},
\EQN
where $\wF \in { \Nhh }, \wF\ge f$ is such that it solves the following minimization problem
\BQN \label{OP}
\cEE{\text{ find the unique $\wF\in { \Nhh }$   so that  } \inf_{ g,f \in { \Nhh }, g \ge f}  
 \norm{g}_{ \Nhh }= \norm{\wF}_{ \Nhh }.}
\EQN
Clearly, \eqref{LD} does not show how to find $\wF$, however it is very helpful for the derivation of upper and lower bounds for $P_f$
since it can be used to check their validity (at least asymptotically), and moreover, it gives further ideas how to proceed. \\
In this paper, for $f\in \kH$ with $f(x_0)>0$ for some $x_0>0$, we find explicitly for $H>1/2$ the unique solution $\wF\in \kH$ of the minimization problem \eqref{OP};
for $H=1/2$ this has already been done in \cite{MR2016767}.
\rE{For the case  $H \in (0,1/2)$,    we determine again $\wF$ under the assumption that $\wF> f$}.
By making use of the Girsanov formula for fBm, 
 we derive in the main result presented in Theorem \ref{ThA} upper and lower bounds for $P_f$.

The paper is organized as follows: Section 2 briefly reviews some results from fractional calculus and related Hilbert spaces. We   introduce weighted fractional integral operators,  fractional kernels and briefly discuss the corresponding reproducing kernel Hilbert spaces.
The main result is presented in Section 3.    Specific properties of fBm that are used in the proof of the main result are displayed in Section 4 followed then by two examples of the drifts, for $H>1/2$ and for $H<1/2$, when the main result holds.
Proofs are   relegated to Section 5. A short Appendix concludes the article.

\section{Preliminaries}
This section reviews basic Riemann-Liouville fractional calculus; a classical reference on this topic is \cite{SamkoKM}.
 We use also the notation and results from \cite{nualart},
 \cite{bioks}, and \cite{Jost}.  We proceed then with the RKHS of fBm.

\begin{defn}
Let $\alpha>0, T>0$. The (left-sided) Riemann-Liouville fractional integral operator of order $\alpha$ over interval $[0,T]$ (or over $\R_+$) is defined by
$$\left(I_{0+}^{\alpha}f\right)(t)=\frac{1}{\Gamma(\alpha)}\int_0^t (t-u)^{\alpha-1}f(u)du,\quad t\in[0,T] \quad (t\in \R_+), $$
where $\Gamma(\cdot)$ is the Euler gamma function. The corresponding right-sided integral operator on $[0,T]$ is defined by
$$\left(I_{T-}^{\alpha}f\right)(t)=\frac{1}{\Gamma(\alpha)}\int_t^T (u-t)^{\alpha-1}f(u)du, \quad t\in[0,T],   
$$
and the right-sided integral operator  
 on $\R_+$ (also known as the Weyl fractional integral operator) 
 is defined by
$$\left(I_{\infty-}^{\alpha}f\right)(t)=\frac{1}{\Gamma(\alpha)}\int_t^\infty (u-t)^{\alpha-1}f(u)du,\quad t\in \R_+. $$
\end{defn}

Throughout the paper,  
 we suppose that $(I_{T-}^{\alpha}f)(t)=0,$ for $t>T$. Note that in the case $u^\alpha f(u)\in L_1(\R_+)$, 
  the integral $(I_{\infty-}^{\alpha}f)$   
  exists and belongs to $L_1(\R_+)$.

  Next, for $p\geq1$, denote
\begin{equation*}
I_+^{\alpha}(L_p[0,T])=\{f: f=I_{0+}^{\alpha}\varphi \text{ for some } \varphi\in L_p[0,T]\},
\end{equation*}
\begin{equation*}\label{Eqn:2}
I_-^{\alpha}(L_p[0,T])=\{f: f=I_{T-}^{\alpha}\varphi \text{ for some } \varphi\in L_p[0,T]\},
\end{equation*}
and define similarly $I_-^{\alpha}(L_p(\R_+))$.
%
%
\cE{If} $0<\alpha<1$, then the function $\varphi$ used in the above definitions 
(it is determined uniquely) coincides for almost all (a.a.)  
 $t\in[0,T]$ $(t\in \R)$ with the left- (right-) sided Riemann-Liouville fractional derivative of $f$ of order $\alpha$. The derivatives are denoted by

$$(I_{0+}^{-\alpha}f)(t)=(\mathcal{D}_{0+}^{\alpha}f)(t)=\frac{1}{\Gamma(1-\alpha)}\frac{d}{dt}\left( \int_0^t(t-u)^{-\alpha}f(u) du \right),$$
$$(I_{\infty-}^{-\alpha}f)(t)=(\mathcal{D}_{\infty-}^{\alpha}f)(t)=-\frac{1}{\Gamma(1-\alpha)}\frac{d}{dt}\left( \int_t^{\infty}(u-t)^{-\alpha}f(u) du \right),$$  
and
$$I_{T-}^{-\alpha}(t)=(\mathcal{D}_{T-}^{\alpha}f)(t)=(\mathcal{D}_{\infty-}^{\alpha}f1_{[0,T]})(t),$$
respectively. Let $f\in I_-^{\alpha}(L_p(\R))$ or $I_{\pm}^{\alpha}(L_p[0,T]), p\geq1, 0<\alpha<1$. Then for the corresponding indices $0,T$, and $\infty$,  
 we have
$$I_{\pm}^{\alpha}\mathcal{D}_{\pm}^{\alpha}f=f.$$ In the case when  $f\in L_1(\R_+)$,   
  we have $\mathcal{D}_{\pm}^{\alpha}I_{\pm}^{\alpha}f=f$ (\cite{SamkoKM})
In the following we introduce weighted fractional integral operators,  fractional kernels and briefly discuss the corresponding reproducing kernel Hilbert spaces.
\def\CAh{\rE{C_1}}
\def\CAhI{\bE{C_1^{-1}}}
Introduce weighted fractional integral operators by

$$(K_{0+}^H f)(t)= \CAh t^{H-1/2}(I_{0+}^{H-1/2}u^{1/2-H}f(u))(t),$$
$$(K_{0+}^{H,*} f)(t)= \CAh^{-1} t^{H-1/2}(I_{0+}^{1/2-H}u^{1/2-H}f(u))(t),$$
$$(K_{\infty-}^H f)(t)= \CAh t^{1/2-H}(I_{\infty-}^{H-1/2}u^{H-1/2}f(u))(t),$$  
and
$$(K_{\infty-}^{H,*} f)(t)= \CAh^{-1} t^{1/2-H}(I_{\infty-}^{1/2- H}u^{H-1/2}f(u))(t),$$
where $ \CAh =\left(\frac{2H\Gamma(H+1/2)\Gamma(3/2-H)}{\Gamma(2-2H)}\right)^{1/2}$.
For $H=\frac{1}{2}$ we put $K_{0+}^H=K_{0+}^{H,*}=K_{\infty-}^{1/2}=K_{\infty-}^{1/2,*}=\I$, where $\I$ is the identity operator.

Let $H>\frac12$. If $u^{H-\frac12}f(u)\in L_1(\R_+)$,  
   then   $K_{\infty-}^{H,*}K_{\infty-}^{H}f=f$. Furthermore, for   $H<\frac{1}{2} $ and for such $f$  that $u^{\frac{1}{2}-H}f(u)\in L_1(\R_+) $,  
   we have that $K_{\infty-}^{H}K_{\infty-}^{H,*}f=f$.  For ${H>\frac{1}{2}}$ and for such $f$ that $u^{H-\frac{1}{2}}f(u)\in I_{-}^{H-\frac{1}{2}}(L_p(\R_+)) $ for some $p\geq 1,$ we have that $K_{\infty-}^{H}K_{\infty-}^{H,*}f=f$. For $f\in L_2(\R_+)$ and $H\in(0,1)$,  
    $K_{0+}^{H}K_{0+}^{H,*}f=f$.

Denote $K_T^H f=K_{\infty-}^H(f 1_{[0,T]})$ and $K_T^{H,*}f=K_{\infty-}^{H,*}(f 1_{[0,T]})$. For $H\in(0,1)$ and $t > s$, define the 
fractional kernel
$$K_H(t, s) := \rE{\frac{ \CAh }{\Gamma\Big(H+1/2\Big)}}\Big(\Big(\frac ts\Big)^{H-\frac12}(t-s)^{H-\frac12}-(H-\frac12)s^{\frac12-H}\int_s^t(u - s)^{H-\frac12}u^{H-\frac32}du\Big).  
$$
 For $H>\frac12$, the  kernel $K_H(t, s)$   
 is simplified to $$K_H(t, s)
 =\frac{ \CAh }{\Gamma\Big(H-\frac12\Big)}s^{\frac12-H}\int_s^t(u - s)^{H-\frac32}u^{H-\frac12}du.$$

In turn, introduce the fractional kernel
$$K_H^*(t, s) = \frac{1}{\rE{\CAh \Gamma(H+1/2)}}\Big(\Big(\frac ts\Big)^{H-\frac12}(t-s)^{\frac12-H}-(H-\frac12)s^{\frac12-H}\int_s^t(u - s)^{\frac12-H}u^{H-\frac32}du\Big).
$$
For $H<\frac12$, the  kernel $K_H^*(t, s)$  
 is simplified to $$K_H^*(t, s)= \frac{s^{\frac12-H}}{\rE{ \CAh \Gamma(1/2-H)} }\int_s^t(u - s)^{-H-\frac12}u^{H-\frac12}du .$$
 \rE{By direct calculations we obtain}
 \begin{equation*}
 (K_{\infty-}^H 1_{[0,t]})(s)=(K_t^H 1_{[0,t]})(s)=K_H(t, s)
  \end{equation*}

and
\begin{equation*}
(K_{\infty-}^{H,*} 1_{[0,t]})(s)=(K_{t}^{H,*} 1_{[0,t]})(s)=K_H^*(t, s).
\end{equation*}
From the   integration-by-parts formula for fractional integrals
$$\int_a^b g(x)I_{a+}^\alpha f(x)dx=\int_a^b f(x)I_{b-}^\alpha g(x)dx$$
for $f\in L_p[a,b]$, $g\in L_q[a,b]$ with $\frac1p+\frac1q\leq 1+\alpha$, we get that for $H>\frac12$ and $f\in L_p[0,t]$ with $p>1$
\begin{equation*}\int_0^t (K_{\infty-}^{H} 1_{[0,t]})(s)f(s)ds=\int_0^t (K_{0+}^{H} f)(s)ds
\end{equation*}
 and for $H<\frac12$ and $f\in L_p[0,t]$ with $p>1$
\begin{equation*}
\int_0^t (K_{\infty-}^{H,*} 1_{[0,t]})(s)f(s)ds=\int_0^t (K_{0+}^{H,*} f)(s)ds.
\end{equation*}
\rE{Next,} we introduce the RKHS of fractional Brownian motion (corresponding results for finite interval are
 described in detail in \cite{decreusefond},
 \cite{nualart}, and \cite{bioks}).  
 Let \rE{$H \in (0,1)$ be fixed and}  {recall that $R_H$ defined in (\ref{Rh}) can be defined also as follows}  
$$R_H(t, s) =\int_0^{t\wedge s}K_H(t, u)K_H(s, u) du.$$
{ \begin{defn}
 (\cite{bioks}) The reproducing kernel Hilbert space (RKHS) of the fractional Brownian motion \rE{on $[0,T]$}, denoted
by $\kH[0,T]$ is defined as the closure
of the vector space spanned by the set of functions ${R_H(t, \cdot), t\in [0,T]}$ with
respect to the scalar product
$\langle R_H(t, \cdot),R_H(s, \cdot)\rangle = R_H(t, s)$, $t, s \in[0,T].$
 \end{defn}
\rE{In \cite{decreusefond} it is shown that} 
  $\kH[0,T]$ is the set of functions $f$ which
can be written as
$f(t) =\int_0^tK_H(t, s)  {\phi} (s) ds$
for some $  {\phi} \in L_2([0, T]).$ By definition,
$\|f\|_{\kH[0,T]}=\|  {\phi} \|_{L_2[0, T]}.$
Extending this definition to $\R_+$, we get the following definition of the RKHS $\kH:=\kH_H(\R_+)$.}
For any $H\in(0,1)$,  
 $ \Nhh $ is the set of functions $f$ which
can be written as
\BQN
f(t) =\int_0^t K_H(t, s)  {\phi}  (s) ds=\int_0^t(K_{\infty-}^H 1_{[0,t]})(s)  {\phi} (s) ds
=\int_0^t (K_{0+}^{H}   {\phi} )(s)ds
\EQN
for some $\phi \in L_2(\R_+).$    
Note that $f'(t)=(K_{0+}^{H}   {\phi} )(t)$ \rE{and} $  {\phi} (t)=(K_{0+}^{H,*}  {f'})(t)$,
 therefore
 $$\|f\|_{ \Nhh }=\|  {\phi} \|_{L_2(\R_+)}=\|K_{0+}^{H,*}  {f'}\|_{L_2(\R_+)}.$$
Next,
define the spaces $L_2^H(\R_+)$ in the following way: 
$$L_2^H(\R_+)=\{f: K_{\infty-}^H|f|\in L_2(\R_+)\}.$$

If  $H\in (0, \frac{1}{2})$, 
we define further
$$\widetilde{L}_2^H(\R_+)=L_2^H(\R_+)\cap \Biggl\{f:\R_+\to\R:\int_0^T t^{1-2H}
\left(\int_T^{\infty}u^{H-1/2}f(u)(u-t)^{H-3/2}du \right)^2 dt \rightarrow 0 \text{ as } T\rightarrow\infty \Biggr\}.$$
For function $g$ that admits the representation $g(t)=\int_0^t g'(s)ds$ introduce the norm
\BQN\norm{g}=\norm{g'}_{L_2(\R_+)}.
\EQN

\newcommand{\Abs}[1]{\Bigl\lvert #1 \Bigr\rvert}

\section{Main result}
In this section we study the boundary non-crossing probability
$P_f=\pk{B^H(t) + f(t) \le u(t),t\in \R_+}$ for $f\in  \HhH $ and a measurable function $u:\R_+ \to \R$  with $u(0)\ge 0$. Throughout this paper,  
 we assume that
$P_0=\pk{B^H(t)\le u(t), t\in \R_+} \in (0,1)$. In applications, see, e.g., \cite{1103.60040,1079.62047}  
 it is of interest to calculate the rate of decrease to 0 of $P_{\gamma f}$ as  $\gamma \to \IF$
for some $f\in  \HhH $. On the other side,  if $\norm{f}_{ \Nhh }$ is small, we expect that $P_f$ is close to $P_0$. Set below
$\alpha= \Phi^{-1}(P_0)$ where $\Phi$ is the distribution function of a $N(0,1)$ random variable. Our first result derives upper and lower bounds of $P_f$
for any $f\in  \HhH $.\\
\def\aF{\bE{g}}
\begin{lem} \label{LemUL} For any $f\in  \HhH $
we have
 \BQN\label{eq:00:2b}
\Abs{P_f - P_0} &\le \frac {1 }{\sqrt{2 \pi}} \norm{f}_{ \Nhh }.
\EQN
If further $\aF\in  \HhH $ is such that $\aF \ge f$, then
\BQN\label{eq:WL}
\Phi(\alpha - \norm{\aF}_{ \Nhh }) \le P_{\aF}\le P_f \le \Phi(\alpha+ \norm{f}_{ \Nhh }).
\EQN
\end{lem}
Clearly, \eqref{eq:00:2b} is useful only if $\norm{f}_{ \Nhh }$ is small. On the contrary, \rE{the lower bound of} \eqref{eq:WL} is important
for $f$ such that $\norm{f}_{ \Nhh }$ is large and
 $\norm{\aF}_{ \Nhh }>0$. 
{Taking $g=\widehat f$, with  $\widehat{f}$ being the solution of (3)}  {and noting that for any $\gamma >0$ we have $\widehat{\gamma f}= \gamma \widehat{f}$ for any $f\in \kH$, then the lower bound in \eqref{eq:WL} implies the following result:}
\begin{cor} \label{corLDP} For any $f\in  \HhH $ such that $f(x_0)>0$ for some $x_0\in (0,\IF)$ the claims in \eqref{LD}
 and \eqref{OP} hold.
\end{cor}

The main result of this section is Theorem \ref{ThA} below which presents upper and lower bounds for $P_f$ under some restriction on $f$ and a general  
measurable $u$ as above.  Let the function $f$ {be} differentiable with  derivative $f'\in L_2(\R_+)$. Then  the operator $(K_{0+}^{H,*} f')$  is \rE{well-defined}. Consider the following assumptions on $f$:

\begin{itemize}
\item[(i)] $(K_{0+}^{H,*} f')\in L_2(\R_+)$, i.e., $f\in \HhH .$

\item[(ii)] Let $h(t):=  \int_0^t (K_{0+}^{H,*} f')(s)ds$. We assume that the smallest concave nondecreasing majorant $\widetilde{h}$ of the function $h$ has the right-hand derivative $\widetilde{h}'$ such that  $\widetilde{h}'\in L_2(\R_+)$ and moreover the function $$K(t):=(K^{H,*}_{\infty-}\widetilde{h}')(t)$$ is nonincreasing, $K\in L_2^H(\R_+)$ for $H>\frac{1}{2}$ and $K\in \widetilde{L}_2^H(\R_+)$ for $H<\frac{1}{2}$, $$K(t)=o(t^{-H})\;\text{as}\; t\to\infty.$$

\item[(iii)] The function $\widetilde{h}'$ can be presented as $\widetilde{h}'(t)=(K_{0+}^{H,*}  \hat{f}' )(t),\;t\in\R_+$,  
 for some $  \hat{f}' \in L_2(\R_+)$. Evidently, in this case the function $\widetilde{h}$
  admits the representation $\widetilde{h}(t) =  \int_0^t (K_{0+}^{H,*}  \hat{f}'  )(s)ds$.  
 Denote $\widehat{f}(t)=\int_0^t \hat{f}' (s)ds=\int_0^t(K_{0+}^{H} \widetilde{h}')(s)ds$.

\end{itemize}

\begin{thm} \label{ThA} 1. Under assumptions $(i)$--$(iii)$  {we have $\widehat{f}\in  \HhH  $ and}
\BQN\label{thA:1}
P_f \le P_{f- \wwF }\exp\left(\int_0^\infty u(s)d(-K(s))-\frac{1}{2}\norm{\widetilde{h}}^2 \right).
\EQN
\cE{2. Suppose  that  $u_{-}: \R_+ \to \R$ is such that $u_{-}(t)< u(t), t\in \R_+$.
 If $H< 1/2$,  
  assume additionally that $ \wwF  \ge f$.  
    Then for any $H\in (0,1)\setminus \{1/2\}$,  
\BQN\label{thA:2}
P_f \ge P_{ \wwF } \geq P( u_{-}(t) \le B^H(t)\leq u(t),t\in\R_+)
\exp\left(\int_0^\infty u_{-}(s)d(-K(s))-\frac{1}{2}\norm{\widetilde{h}}^2  \right)
\EQN
holds, provided that $\int_0^\infty u_{-}(s)d(-K(s))$ is finite.}
\end{thm}

As we show below,  
 the upper and lower bounds above become (in the log scale) precise when  ${f}$ is large.
\cEE{\begin{cor} \label{crA} Under the assumptions
and notation of Theorem \ref{ThA}, if further $f(x_0)>0$ for some $x_0 \in (0,\IF)$, then
\BQN\label{crA:1}
- \ln P_{\gamma f}  \sim  \frac{\gamma^2}{2} \norm{\widetilde{h}}^2, \quad \gamma \to \IF.
\EQN
\end{cor}
}

As a by-product, we solve the minimization problem \eqref{OP}, namely we have
\cEE{\begin{cor} \label{crB} Under the assumptions
and notation of Theorem \ref{ThA}
\BQN\label{crB:1}
\cEE{\inf_{f,g\in { \Nhh },g \ge f} \norm{g}_{ \Nhh } =  \norm{\wwF}_{ \Nhh }= \norm{\widetilde{h}}.}    
\EQN
\end{cor}
}

{\bf Remarks}: a) If $H\in (1/2, 1)$, then under conditions $(i)$--$(iii)$,  
 we find that $\wwF$ is the explicit solution of the minimization problem \eqref{OP}.\\
b) The case $H=1/2$ is discussed in \cite{BiHa1}, see also \cite{MR2016767}.\\
c) It follows from Lemma \ref{AppE} that  for $H>\frac12$,
 $ \wwF  \ge f$ because it immediately follows from this lemma and inequality $\tilde{h}\geq h$ that $\hat{f}'\geq f'$.

\section{Auxiliary  results}
For the proof of our main result,   
 we need to discuss several properties of fBm.  This section discusses first
the relation between fBm, Molchan martingale and the underlying Wiener process. Then we consider the Girsanov theorem which is crucial for our analysis.

 \subsection{Fractional Brownian motion, Molchan martingale and ``underlying'' Wiener process}

In what follows we consider  continuous modification of fBm that exists due to  well-known Kolmogorov's theorem.  
Denote by $\F^{B^H}=\{\F_t^{B^H},t\in \R_+\}$ with $\F_t^{B^H}=\sigma\{B^H(s),0\leq s\leq t\}$ the filtration generated by $B^H$.
\rE{Below} we establish the following relation. According to  \cite{nualart}, \cite{bioks}, \cite{Jost}, and \cite{NVV}, $B^H$ can be presented as  

\begin{equation}\label{Eqn:4}
B^H(t)=\int_0^t(K_{\infty-}^H 1_{[0,t]})(s)dW(s)=\int_0^t(K_t^H 1_{[0,t]})(s)dW(s)=\int_0^t K_H(t,s)dW(s),
\end{equation}
where $W=\{W(t),t\in \R_+\}$ is an ``underlying'' Wiener process whose filtration coincides with $\F^{B^H}$.  Evidently, \begin{equation}\label{Eqn:5}
W_t=\int_0^t(K_{\infty-}^{H,*}1_{[0,t]})(s)dB^H(s)=\int_0^t(K_t^{H,*}1_{[0,t]})(s)dB^H(s)=\int_0^tK_{H}^*(t,s)dB^H(s).
\end{equation}
 Another form of relations (\ref{Eqn:4}) and (\ref{Eqn:5}) can be obtained in the following way. According to \cite{NVV}, we can introduce the kernel
\BQN\label{LH}
l_H (t,s)=\rE{\left(\frac{\Gamma(3-2H)}{2H\Gamma(3/2-H)^3\Gamma(H+1/2)}\right)^{1/2}}s^{1/2-H}(t-s)^{1/2-H} 1_{[0,t]}(s) ,\quad s, t\in\R_+
\EQN
and consider the process
\begin{equation}\label{Eqn:7}
M^H(t)=\int_0^t l_H(t,s)dB^H(s),\quad t \in\R_+,\quad H\in(0,1).   
\end{equation}

The process $M^H$ from (\ref{Eqn:7}) defines a Gaussian square-integrable martingale with square characteristics $\langle M^H\rangle(t)=t^{2-2H}$, $t\in\R_+$,   
 and with filtration $\F^{M^H}\equiv\F^H$.  Then    
  the process $\widetilde{W}(t)=(2-2H)^{-1/2}\int_0^t s^\alpha dM^H(s)$ is a Wiener process with the same filtration.

\begin{lem} \label{Lem1} The processes $\widetilde{W}$ and $W$ coincide.
\end{lem}

\begin{defn} (\cite{bioks}, \cite{Jost}, \cite{nualart}) Wiener integral w.r.t. fBm is defined  for any $T\in \R_+$ and $H\in(0,1)$ as
\BQNY
\int_0^Tf(s)dB^H(s)&=&\int_0^T(K_{\infty-}^Hf 1_{[0,T]})(s)dW(s)                            
=\int_0^{\infty}(K_{\infty-}^Hf 1_{[0,T]})(s)dW(s)\\
&=&\int_0^{\infty}(K_T^Hf)(s)dW(s)=\int_0^T(K_T^Hf)(s)dW(s)
\EQNY
and the integral $\int_0^Tf(s)dB^H(s)$ exists for $f\in L_2^H (\R_+).$
\end{defn}

Now we extend the notion of integration w.r.t.\ fBm on the $\R_+$ from $[0,T]$ by the following definition.

\begin{defn}\begin{equation}\label{Eqn:8}\int_0^\infty f(s)dB^H(s)=L_2\text{-}\Lim_{T\to\infty}\int_0^T f(s)dB^H(s),\end{equation}
if this limit exists.
\end{defn}

\begin{lem}\label{Lemma 2.2} Let function $f\in\L_2^H(\R_+)$ for $H>\frac{1}{2}$ and $f\in\widetilde{L}_2^H(\R_+)$ for $H<\frac{1}{2}$.

Then the limit in the right-hand side of (\ref{Eqn:8}) exists \bE{and} 
\BQN
\int_0^\infty f(s)dB^H(s)=\int_0^\infty (K^H_{\infty-} f)(s)dW(s).
\EQN
\end{lem}


\begin{lem}\label{Lemma2.33} Let $h=h(t),t\in R_+$,    
be a nonrandom measurable function \bE{such that} 
\begin{enumerate} \item $h\in L_2^H(\R_+)$ for $H>\frac{1}{2}$ and $h\in \widetilde{L}_2^H(\R_+)$ for $H<\frac{1}{2}$; \item $h$ is nonincreasing; \item $s^H h(s)\to 0 $ as $s\to\infty.$
\end{enumerate}
Then there exists integral $\int_0^\infty h(s)dB_s^H$ in the sense of Lemma 2.2 and \bE{moreover}
\BQN
\int_0^\infty h(s)dB^H(s)=\int_0^\infty B^H(s)d(-h(s)),
\EQN
 where the integral in the right-hand side is a Riemann-Stieltjes integral with continuous integrand and nondecreasing integrator.
\end{lem}
\COM{{Next, we introduce few elements from convex analysis.}
Let $h:\R_+\to\R$ be a measurable function from the space $\mathcal{H}_{\frac{1}{2}}.$ Denote $\widetilde{h}$ the smallest nondecreasing concave majorant of $h$.   
According to \cite{BiHa1}, $\widetilde{h}\in\mathcal{H}_{\frac{1}{2}},$ $\widetilde{h}\geq0$ and $\widetilde{h}$ for a.a.\   
 $t\in\R_+$ has a derivative $\widetilde{h}'$.
Furthermore, we can assume that $\widetilde{h}'$ is the right-hand derivative. As it was mentioned in \cite{BiHa1}, $\Inf_{g\geq h,g\in\mathcal{H}_{\frac{1}{2}}}\norm{g}=\norm{\widetilde{h}}$ and $\norm{h}^2=\norm{\widetilde{h}}^2+\norm{h-\widetilde{h}}^2.$
}
\subsection{Girsanov theorem for fBm} Let $H\in(0,1)$. Consider a fBm with absolutely continuous  drift $f$ that admits a following representation: $ B^H(t)+f(t)=B^H(t)+\int_0^t f'(s)ds.$
To annihilate the drift, there are two  equivalent approaches. The first one  is to assume that $K_{H}^*(t,\cdot)f'(\cdot)=(K_{0+}^{H,*}f')(\cdot)\in L_1[0,t]$ for any $t\in R_+$, to  equate \begin{equation}\label{Eqn:11}B^H(t)+f(t)=\widehat{B}^H(t),\end{equation}
where $\widehat{B}^H$ is the fBm with respect to the new probability measure,
and accordingly to \eqref{Eqn:5},     
 to transform (\ref{Eqn:11}) as
$$\int_0^t(K_{\infty-}^{H,*}1_{[0,t]})(s)dB^H(s)+\int_0^t(K_{\infty-}^{H,*}1_{[0,t]})f'(s)ds
=\int_0^t(K_{\infty-}^{H,*}1_{[0,t]})(s)d\widehat{B}^H(s),$$
or,
$$\int_0^tK_{H}^*(t,s)dB^H(s)+\int_0^tK_{H}^*(t,s)f'(s)ds
=\int_0^tK_{H}^*(t,s)d\widehat{B}^H(s),$$
or, at last,
$$W(t)+\int_0^t(K_{\infty-}^{H,*}1_{[0,t]})(s)f'(s)ds=W(t)+\int_0^t(K_{0+}^{H,*}f')(s)ds=\widehat{W}(t),$$
where $\widehat{W}=\{\widehat{W}_t,t\in\R_+\}$ is a Wiener process with respect to a new probability measure $Q$, say. The second one is  to apply  Girsanov's 
theorem from \cite{Yu08LNotes}. We start with (\ref{Eqn:11}); suppose that $s^{\frac12-H}f'(s)\in L_1[0,t]$ for any $t\in\R_+$ and transform (\ref{Eqn:11}) as follows  {(recall $l_H$ is defined in \eqref{LH})}:
$$M^H(t)+\int_0^t l_H(t,s)f'(s)ds=\widehat{M}^H(t).$$

Further, suppose that the function $q(t)=\int_0^t l_H(t,s)f'(s)ds$ admits the representation
\begin{equation}\label{equ.2.23}
q(t)=\int_0^t q'(s)ds.
\end{equation}  Then
$$(2-2H)^{\frac{1}{2}}\int_0^t s^{\frac{1}{2}-H}dW(s)+\int_0^t q'(s)ds=(2-2H)^{\frac{1}{2}}\int_0^t s^{\frac{1}{2}-H}d\widehat{W}(s),$$
whence $ W(t)+{(2-2H)^{-1/2}}\int_0^tq'(s)s^{H-\frac{1}{2}}ds=\widehat{W}(t).$ Evidently, if the representation \eqref{equ.2.23} holds, then \begin{equation}\label{equ.2.24}{(2-2H)^{-1/2}}
\int_0^tq'(s)s^{H-\frac{1}{2}}ds=\int_0^tK_{H}^*(t,s)f'(s)ds=\int_0^t(K_{0+}^{H,*}f')(s)ds.\end{equation}
Now we give simple sufficient conditions  of existence of $q'$ and $\int_0^tq'(s)s^{H-\frac{1}{2}}ds$.  {The proof} consists in differentiation and integration by parts therefore \bE{it} is omitted.
\begin{lem} \label{lemdrif} (i) Let $H<\frac12 $. Suppose that {the} drift $f$ is absolutely continuous and for any $t>0$,   the derivative 
$|f'(s)|\leq C(t)s^{H-\frac32+\varepsilon}$, $s\leq t$, for some $\varepsilon>0$ and some nondecreasing function $C(t):\R_+\rightarrow\R_+$.
Then for any $t>0$,   
 $$q'(t)={\left(\frac{\Gamma(3-2H)}{2H\Gamma(3/2-H)^3\Gamma(H+1/2)}\right)^{1/2}}\int_0^t s^{1/2-H}(t-s)^{-1/2-H}f'(s)ds $$ and \eqref{equ.2.24} holds.

  (ii) Let $H>\frac12 $. Suppose that the drift $f$   
  is absolutely continuous. Also, suppose that
   there exists the continuous  derivative $(s^{\frac12-H}f'(s))'$ and  
    $ (s^{\frac12-H}f'(s))'\rightarrow 0$ as  $s\rightarrow 0$.
    Then for any $t>0$
    $$q'(t)={\left(\frac{\Gamma(3-2H)}{2H\Gamma(3/2-H)^3\Gamma(H+1/2)}\right)^{1/2}}\int_0^t (t-s)^{1/2-H}(s^{\frac12-H}f'(s))'ds $$ and \eqref{equ.2.24} holds.
\end{lem}

{For a drift $f$ as in Lemma \ref{lemdrif}, then} 
$B^H(t)+\int_0^t f'(s)ds$ is fBm $\widehat{B}^H(t),$ $t\in\R$, say,  
 under such measure $Q$ that
\begin{equation}\begin{gathered}
\label{Eqn:12}\frac{dQ}{dP}=\exp\Big(-\int_0^\infty (K_{0+}^{H,*}f')(s)dW(s)-\frac{1}{2}\int_0^\infty|(K_{0+}^{H,*}f')(s)|^2 ds\Big)\\=\exp\Big(-\int_0^\infty (K_{0+}^{H,*}f')(s)dW(s)-\frac{1}{2}\|f\|^2_{ \Nhh }\Big)\end{gathered}\end{equation}
if (\ref{Eqn:12}) defines a new probability measure.
So, we get the following result.

\begin{thm}\label{thm4.1} If $f\in  \Nhh $, then  $B^H(t)+\int_0^t f'(s)ds=\widehat{B}^H(t), $ where $\widehat{B}^H(t)$ is a fBm under \bE{a} measure $Q$ that satisfies relation \eqref{Eqn:12}.
\end{thm}

\section{Examples of admissible drifts}

We present next two examples of drifts satisfying conditions $(i)$-$(iii)$.

\begin{exm}\label{exa4.1} In order to construct the drift, we start with $h$ and $\widetilde{h}$. Let $H>\frac{1}{2}$, $h(t)=\widetilde{h}(t)=\int_0^ t s^{1/2-H}e^{-s}ds.$ Note that $\widetilde{h}'(t)=t^{1/2-H}e^{-t},$ $t>0$, $\widetilde{h}'\in L_2(\R_+),$ $\widetilde{h}'>0$ and decreases on $\R_+$, therefore $\widetilde{h}$ is a concave function as well as $h$, and evidently, $\widetilde{h}$ is a smallest concave nondecreasing majorant of $h$. Further,
\BQNY
(K_{\infty-}^{H,*}\widetilde{h}')(t)&=&- \CAhI t^{1/2-H}\frac{d}{dt}\left(\int_t^\infty(\bE{z}-t)^{1/2-H}e^{-\bE{z}}d\bE{z}\right)\\
&=&- \CAhI t^{1/2-H} \frac{d}{dt}\left(\int_0^\infty {\bE{z}}^{1/2-H}e^{-\bE{z}-t}d\bE{z}\right)\\
&=& \CAhI \Gamma\Big(\frac{3}{2}-H\Big)t^{\frac{1}{2}-H}e^{-t} = \CAhI \Gamma\Big(\frac{3}{2}-H\Big)\widetilde{h}'(t).
\EQNY
Consequently, the function $K(t):=(K_{\infty-}^{H,*}\widetilde{h}')(t)$ is nonincreasing, $$K_{\infty-}^H(K_{\infty-}^{H,*}\widetilde{h}')(t)=\widetilde{h}'(t)\in L_2(\R_+)$$ implying thus $K\in L_2^H(\R_+)$ and moreover, $K(t)t^H\to 0$ as $t\to \infty$. It means that condition $(ii)$ holds.

Denote $f$, yet the  unknown drift, and let $q(t)=C_2\int_0^t s^{1/2-H}(t-s)^{1/2-H}f'(s)ds,$ \rE{with $C_2:=\rE{\frac{ \CAh }{\Gamma\Big(H+1/2\Big)}}$.}
Then $s^{H-1/2}q'(s)=h'(s)=s^{1/2-H}e^{-s},$ $q'(s)=s^{1-2H}e^{-s}$ and
\BQNY C_2\int_0^t(t-s)^{1/2-H}s^{1/2-H}f'(s)ds=\int_0^t s^{1-2H}e^{-s}ds
\EQNY
and hence  with $C_3=C_2B(\frac{3}{2}-H,H-\frac{1}{2})$ \rE{where $B(a,b)=\Gamma(a)\Gamma(b)/\Gamma(a+b)$}, we obtain  
$$(H-\frac{1}{2})C_3\int_0^t s^{1/2-H}f'(s)ds=\int_0^t (t-s)^{H-1/2}s^{1-2H}e^{-s}ds$$
implying that
$$f(t)={\left(\frac{\Gamma\Big(\frac32-H\Big)}{2H\Gamma(2-2H)\Gamma\Big(H+\frac12\Big)}\right)^{-\frac12}}
\int_0^t s^{H-\frac{1}{2}}\int_0^s(s-\bE{z})^{H-\frac{3}{2}}{\bE{z}}^{1-2H}e^{-\bE{z}}d\bE{z}ds.$$

{Since} $(K_{0+}^{H,*}f')(t)= \CAh t^{H-1/2}q'(t)=t^{1/2-H}e^{-t}\in L_2(\R_+)$ condition $(i)$ holds. Condition $(iii)$ is clearly satisfied since
we can put $\wwF=f$. \cE{Note in particular that the assumption $\wwF\ge f$ if $H \in (0,1/2)$ also holds.}
\end{exm}

\begin{exm}
Let $H<\frac{1}{2}$ and put $h(t)=\widetilde{h}(t)=\int_0^t s^\gamma e^{-s}ds$ with some $0>\gamma>-\frac{1}{2}$ to have $h'$ and $\widetilde{h}'$ in $L_2(\R_+)$. Then, as before, $\widetilde{h}$ is a smallest nondecreasing concave majorant of $h$. Further, we may write
\BQNY
(K_{\infty-}^{H,*}\widetilde{h}')(t)&=& \CAhI t^{1/2-H}\int_t^\infty(\bE{z}-t)^{-H-1/2}{\bE{z}}^{H-1/2+\gamma}e^{-\bE{z}}d\bE{z}\\
&=&- \CAhI t^{1/2-H+\gamma}\int_1^\infty(\bE{z}-1)^{-H-1/2}{\bE{z}}^{H-1/2+\gamma}e^{-\bE{z}t}d\bE{z}
\EQNY
and $K(t):=(K_{\infty-}^{H,*}\widetilde{h}')(t)$ is nonincreasing for $\frac{1}{2}-H+\gamma\leq0$, or $-\frac{1}{2}<\gamma\leq H-\frac{1}{2}.$ Moreover, for $\gamma=H-\frac{1}{2}$
$$|(K_{\infty-}^{H,*}\widetilde{h}')(t)|\leq  \CAhI t^{1/2-H}e^{-t/2}\int_1^\infty {\bE{z}}^{H-3/2}d\bE{z},$$
$$K_{\infty-}^H(|K_{\infty-}^{H,*}\widetilde{h}'|)(t)=\widetilde{h}'(t)\in L_2(\R_+)$$
implying $K\in L_2^H (\R_+)$ and ${\lim_{t\to \IF}}K(t)t^H= 0$. Consequently, condition $(ii)$ holds. Similarly to Example \ref{exa4.1}, $$s^{H-1/2}q'(s)=h'(s)=s^{H-1/2}e^{-s},\;q'(s)=e^{-s}$$ and $C_2\int_0^t(t-s)^{1/2-H}s^{1/2-H}f'(s)ds=\int_0^t e^{-s}ds,$ whence
\begin{equation}\label{Eqn:21}\Big(\frac{1}{2}-H\Big)C_2\int_0^t (t-s)^{-1/2-H}s^{1/2-H}f'(s)ds=1-e^{-t}.\end{equation}

It follows from (\ref{Eqn:21}) that
$$\Big(\frac{1}{2}-H\Big)C_2B(H+\frac{1}{2},\frac{1}{2}-H)\int_0^t s^{1/2-H}f'(s)ds=\int_0^t(t-s)^{H-1/2}(1-e^{-s})ds.$$

Denote $C_4:=(\frac{1}{2}-H)B(H+\frac{1}{2},\frac{1}{2}-H).$  
Then
$$\int_0^t s^{1/2-H}f'(s)ds=\rE{\frac{1}{C_4}}\int_0^t\frac{(t-s)^{H+1/2}}{H+1/2}e^{-s}ds,$$
and
$$t^{1/2-H}f'(t)=\rE{\frac{1}{C_4}}(H)\int_0^t(t-s)^{H-1/2}e^{-s}ds,$$
whence
$$f'(t)=\rE{\frac{1}{C_4}}t^{H-1/2}\int_0^t(t-s)^{H-1/2}e^{-s}ds.$$
Consequently,   
\BQNY
f(t)&=&\rE{\frac{1}{C_4}}\int_0^t s^{H-1/2}\int_0^s(s- {z})^{H-1/2}e^{- {z}}d {z}ds\\
&=&\rE{\frac{1}{C_4}}\int_0^t e^{- {z}}\int_{ {z}}^t s^{H-1/2}(s- {z})^{H-1/2}dsd\bE{z}.
\EQNY
Clearly,  $(K_{0+}^{H,*}f')(t)= \CAh t^{H-1/2}q'(t)=t^{H-1/2}e^{-t}\in L_2(\R_+),$ and condition $(i)$ holds. Condition $(iii)$ is evident.
\end{exm}

\section{Proofs}
\subsection{Proofs of auxiliary results}
{\bf Proof of Lemma \ref{Lem1}}: It was established in \cite{NVV} that fBm $B^H$ can be ``restored'' from $\widetilde{W}$ by the following formula
$B^H(t)=\int_0^t K_H(t,s)d\widetilde{W}(s)$,
 but it means
$$\widetilde{W}(t)=\int_0^t(K_{\infty-}^{H,*} 1_{[0,t]})(s)dB^H(s)=W(t),$$
hence the proof follows. \qed

{\bf Proof of Lemma \ref{Lemma 2.2}}:
 On one hand, we have that $\int_0^\infty(K^H_{\infty-} f)(s)dW(s)$ exists. On the other hand, we have the equality $\int_0^T f(s)dB^H(s)=\int_0^T(K_{\infty-}^H f 1_{[0,T]})(s)dW(s)$. At last,
$$\E{\Biggl(\int_0^\infty(K^H f)(s)dW(s)-\int_0^T(K_{\infty-}^H f 1_{[0,T]})(s)dW(s) \Biggr)^2}$$
\begin{equation}\label{Eqn:9}=\int_T^\infty((K_{\infty-}^H f)(s))^2 ds+\int_0^T((K_{\infty-}^H f)(s)-K_T^H f)(s))^2 ds.\end{equation}

Since $f\in L_2^H(\R_+),$ we have that $\int_T^\infty((K_{\infty-}^H f)(s))^2 ds\to 0,$ $T\to\infty.$ Further, let $H>\frac{1}{2}.$ Then
\begin{equation}\label{Eqn:10}\int_0^T((K_{\infty-}^H f-K_T^H f)(s))^2 ds=\rE{C_1}\int_0^T s^{1-2H}\left(\int_T^\infty f(t)t^{H-\frac{1}{2}}(t-s)^{H-\frac{3}{2}}dt\right)^2 ds.\end{equation}

Since $|f|\in L_2^H (\R_+)$ together with $f$, we have that for any $s\leq T$
$$\int_T^\infty|f(t)|t^{H-\frac{1}{2}}(t-s)^{H-\frac{3}{2}}dt\to 0\text{ as }T\to\infty$$
and is dominated by $\int_s^\infty |f(t)|t^{H-1/2}(t-s)^{H-3/2}dt.$ Therefore, the right-hand side of (\ref{Eqn:10}) tends to $0$ due to the Lebesgue dominated convergence theorem.
Next, for $0<H<\frac{1}{2}$ and by the definition $\widetilde{L}_2^H(\R_+)$, we have  
\BQNY
\int_0^T((K_{\infty-}^H f)(s)-(K_T^H f)(s))^2 ds&=&C_1^2\int_0^T s^{1-2H}\Big(\frac{d}{ds}\Big(\int_s^\infty u^{H-\frac{1}{2}} f(u)(u-s)^{H-\frac{1}{2}}du\Big)\\
&&-\frac{d}{ds}\Big(\int_s^T u^{H-\frac{1}{2}}f(u)(u-s)^{H-\frac{1}{2}}du\Big)\Big)^2ds\\
&=&C_1^2\int_0^T s^{1-2H}\Big(\int_T^\infty u^{H-\frac{1}{2}} f(u)(u-s)^{H-\frac{3}{2}}du\Big)^2 ds\\
& =& C_1^2\int_0^T s^{1-2H}\Big(\int_T^\infty u^{H-\frac{1}{2}}f(u)(u-s)^{H-\frac{3}{2}}du\Big)^2 ds\to 0
\EQNY
as $T\to\infty$ implying that the right-hand side of (\ref{Eqn:9}) vanishes as $T\to\infty$, hence the claim follows. \QED

{\bf Proof of Lemma \ref{Lemma2.33}}: According to Lemma \ref{Lemma 2.2}, under condition 1) the integral $\int_0^\infty h(s)dB^H(s)$ exists,
\begin{equation}\label{Eqn:16}\int_0^\infty h(s) dB^H(s)=\int_0^\infty(K^H_{\infty-} f)(s)dW(s)=L_2\text{-}\lim_{T\to\infty}\int_0^T f(s)dB^H(s).\end{equation}

Further, it was mentioned in \cite{Jost} that $\int_0^T h(s)dB^H(s)$ is an $L_2$-limit of the corresponding integrals for the elementary functions:
\BQN \label{Eqn:17}
\int_0^T h(s) dB^H(s)&=&L_2\text{-}\lim_{|\pi|\to 0}\sum_{i=1}^{N} h(s_{i-1})(B^H(s_i)-B^H(s_{i-1}))\notag\\
&=&L_2\text{-}\lim_{|\pi|\to 0}(\sum_{i=1}^{N}B^H(s_i)(s_{i-1})-h(s_i))+B^H(T)h(T))\notag\\
&=&\int_0^T B^H(s)d(-h(s))+B^H(T)h(T).
\EQN

{In view of}  (\ref{Eqn:16}), the limit   
in the right-hand side of (\ref{Eqn:17}) exists and due to condition 3),  
 it equals $\int_0^\infty B^H(s)d(-h(s)),$ whence the proof follows. \QED

\subsection{{Proofs} of the main results}
 {\it Proof of Lemma \ref{LemUL}}:
If $f=0$, then $\norm{f}_{ \Nhh }=0$, hence the first claim follows. Assume therefore that $\norm{f}_{ \Nhh }>0$.
In view of \cite{Mandjes07} (see page 47 and 48 therein),
 a standard {fBm}  $B_H(t),t \ge 0$ can be
realized in the separable Banach space
$$E=\biggl\{\omega: \R \to \R, \text{  continuous}, \quad \omega(0)=0, \quad  \lim_{t \to \IF} \frac{\lvert \omega(t) \rvert }{1+ t}=0\biggr\}$$
equipped with the norm $\norm{\omega}_E= \sup_{t\ge 0} \frac{\lvert \omega(t) \rvert}{1+ t}$.
Consequently, Theorem 1' in \cite{LiKuelbs} can be applied, hence
\BQN\label{eq:WL0}
\Phi(\alpha - \norm{f}_{ \Nhh }) \le P_f \le \Phi(\alpha+ \norm{f}_{ \Nhh }).
\EQN
Since for any $\aF \ge f$ we have $P_{\aF}\le P_f$,
 then the claim in \eqref{eq:WL} follows.
Next, in view of \eqref{eq:WL0},
 we have by the mean value theorem (see also Lemma 5 in \cite{Janssen08})
\BQNY
P_f- P_0 &\le  & \Phi(\alpha + \norm{f}_{ \Nhh })- \Phi(\alpha)
=  \norm{f}_{ \Nhh } \Phi'(c)\le \frac{\norm{f}_{ \Nhh }}{\sqrt{2 \pi}}
\EQNY
for some real $c$ and similarly using again \eqref{eq:WL0},   
\BQNY
P_f- P_0 &\ge & \Phi(\alpha - \norm{f}_{ \Nhh })- \Phi(\alpha)
\ge -\frac{\norm{f}_{ \Nhh }}{\sqrt{2 \pi}},
\EQNY
hence the proof is complete.
\QED

\cEE{{\it Proof of Corollary \ref{corLDP}:} In view of \eqref{eq:WL} we have for any $\gamma >0$ and any $g\in { \Nhh }, g\ge f$
\BQNY
 P_{\gamma f} \ge  P_{\gamma g}  \ge \Phi(\alpha - \gamma \norm{g}_{ \Nhh }).
\EQNY
Since  $g(x_0)>0$ follows from $f(x_0)>0$, then $\norm{g}_ { \Nhh }>0$, hence for all $\gamma$ large
$$  \ln P_{\gamma f} \ge - (1+o(1))\frac{\gamma^2}{2} \inf_{ g\in { \Nhh }, g\ge f} \norm{g}_{ \Nhh }^2.$$  
Since the norm is a convex function and the set $A_f:=\{g\in { \Nhh }, g\ge f\}$ is convex, then the minimization problem \eqref{OP}
has a unique solution $\wwF$, and thus the proof is complete. \QED}

{\bf Proof of Theorem \ref{ThA}:} Define the function $h(t)=\int_0^t h'(s) ds$ with $$h'(s)=f_H(s):=(K_{0+}^{H,*}f')(s)$$
 and introduce its smallest concave nondecreasing majorant $\widetilde{h}$.
  As shown in \cite{BiHa1} $\widetilde{h}(t)=\int_0^t \widetilde{h}'(s) ds$  
   and
$$\norm{h}^2:=\int_0^\infty (h'(s))^2ds=\int_0^\infty (f_H(s))^2 ds=\|f\|^2_{  \Nhh }  
=\norm{\widetilde{h}}^2 +\norm{h-\widetilde{h}}^2 .$$
Next, let the probability measure $Q$ be defined  by the relation
\begin{eqnarray}\label{Eqn:18}
\frac{dQ}{dP}&=&\exp\left(-\int_0^\infty f_H(s) dW(s)-\frac{1}{2}\norm{h}^2\right)   
= \exp\left(-\int_0^\infty f_H(s) d \widehat{W}(s)+\frac{1}{2}\norm{h}^2\right),
\end{eqnarray}
where $W$ is the ``underlying'' Wiener process, $d\widehat{W}=dW+f_H(s)ds,$ $\widehat{W}$ is a Wiener process w.r.t. the measure $Q$. Note that (\ref{Eqn:18}) defines {a} probability measure since $ f_H\in L_2(\R_+)$, due to $(i)$ and Theorem \ref{thm4.1}. Then
\begin{eqnarray*}
P_f&=&\mathbb{E}_Q \left\{\1\{B^H(t)+f(t)\leq u(t),t\in\R_+\}\frac{dP}{dQ} \right\}\\
&=&\mathbb{E}_Q \left\{\1\{\widehat{B}^H(t)\leq u(t),t\in\R_+\} \exp\left(\int_0^\infty f_H(s)
 d \widehat{W}(s)-\frac{1}{2}\norm{h}^2\right\} \right\}
\\
&=&\mathbb{E}\left\{ \1\{B^H(t)\leq u(t),t\in\R_+\}\exp\left(\int_0^\infty f_H(s)
 d \widehat{W}(s)-\frac{1}{2}\norm{h}^2\right)\right\} .
\end{eqnarray*}

Furthermore,
\BQNY
\int_0^\infty f_H(s)dW(s)=\int_0^\infty( f_H(s)-\widetilde{h}'(s))dW(s)+\int_0^\infty \widetilde{h}'{(s)}dW(s)\\=\int_0^\infty(h'(s)-\widetilde{h}'(s))dW(s)+\int_0^\infty \widetilde{h}'{(s)}dW(s).
\EQNY
Next setting $K(t):=(K^{H,*}_{\infty-}\widetilde{h}')(t)$,  we have   
 $$\int_0^\infty \widetilde{h}'(s)dW(s)=\int_0^\infty(K^{H,*}_{\infty-}\widetilde{h}')(s)dB^H(s)=
\int_0^\infty  {K(s)} dB^H(s) $$ and both integrals are correctly defined.
Indeed, $\widetilde{h}'\in L_2(\R_+)$ implying  that $\int_0^\infty \widetilde{h}'(s)dW(s)$ exists.
Moreover,   
{in view of $(ii)$},  
$$ {K}\in L_2^H(\R_+)$$ for $H>\frac{1}{2}$  and $ {K}\in \widetilde{L}_2^H(\R_+)$ for $H<\frac{1}{2}$,
 therefore $\int_0^\infty  {K}(s)dB^H(s)$ exists, according to Lemma \ref{Lemma 2.2} and, furthermore,
 equality (\ref{Eqn:17}) holds. {In the light of} Lemma \ref{Lemma2.33},  
  we get $$\int_0^\infty K(s)d  {B^H(s)}=\int_0^\infty B^H(s) d(-K(s)).$$
  Consequently, condition $(iii)$ implies (set $I_{K,u}:= \int_0^\infty u(s)d(-K(s)))$
\BQNY
P_f&=&\mathbb{E}\Biggl\{\1\{B^H(t)\leq u(t),t\in\R_+\}
\exp\Biggl(\int_0^\infty(h'(s)-\widetilde{h}'(s))dW(s) -\frac{1}{2}\norm{h-\widetilde{h}}^2
+ {\int_0^\infty B^H(s)d(-K(s))}-\frac{1}{2}\norm{\widetilde{h}}^2\Biggr)\Biggr\}\\
 &\cE{\le}&\mathbb{E}\Biggl\{\1\{B^H(t)\leq u(t),t\in\R_+\}
\exp\Biggl(\int_0^\infty(h'(s)-\widetilde{h}'(s))dW(s) -\frac{1}{2}\norm{h-\widetilde{h}}^2 +I_{K,u}-\frac{1}{2}\norm{\widetilde{h}}^2\Biggl)\Biggr\}\\
&=&\exp\Bigg(I_{K,u}-\frac{1}{2}\norm{\widetilde{h}}^2\Bigg) \mathbb{E}\Bigg(\1\{B^H(t)\leq u(t),t\in\R_+\}\\
&& \times \exp\Bigg(\int_0^\infty ((K_{0+}^{H,*}f')(s)-(K_{0+}^{H,*}\widehat{h})(s))dW(s) -\frac{1}{2}\int_0^\infty ((K_{0+}^{H,*}f')(s)-(K_{0+}^{H,*}\widehat{h}))^2ds\Bigg)\Bigg)\\
&=&\exp \Bigg(I_{K,u}-\frac{1}{2}\norm{\widetilde{h}})^2\Bigg)P_{f-\wwF}.
\EQNY
So, the upper bound \eqref{thA:1} is proved. In order to prove \eqref{thA:2}, \cE{note that in view of Lemma \ref{AppE} for $H\in (1/2, 1)$
\BQNY
 \widehat f \ge f,
\EQNY
which is also assumed to hold if $H \in (0,1/2)$. Clearly the above inequality   
implies that $P_f \ge P_{\widehat f}$.
As above, we have for some function $u_{-}(t) < u(t),t\in \R_+$
\BQNY
P_{\widehat f}&=& \mathbb{E}\Bigg\{\1\{B^H(t)\leq u(t),t\in\R_+\}
\exp\Bigg(\cE{\int_0^\infty B^H(s)d(-K(s))}-\frac{1}{2}\norm{\widetilde{h}}^2\Bigg)\Bigg\}\\
&\ge & \mathbb{E}\Bigg\{\1\{u_{-}(t) \le B^H(t)\leq u(t),t\in\R_+\}
\exp\Bigg(\cE{\int_0^\infty B^H(s)d(-K(s))}-\frac{1}{2}\norm{\widetilde{h}}^2\Bigg)\Bigg\}\\
&\ge & \mathbb{E}\Bigg\{\1\{u_{-}(t) \le B^H(t)\leq u(t),t\in\R_+\}
\exp\Bigg(\cE{\int_0^\infty u_{-}(s)d(-K(s))}-\frac{1}{2}\norm{\widetilde{h}}^2\Bigg)\Bigg\}\\
&\ge & P( u_{-}(t) \le B^H(t)\leq u(t),t\in\R_+)
\exp\Bigg(\cE{\int_0^\infty u_{-}(s)d(-K(s))}-\frac{1}{2}\norm{\widetilde{h}}^2\Bigg)\Bigg\},
\EQNY
hence the proof is complete.} \QED

 {{\bf Proof of Corollary \ref{crA}:} \cE{Since $\wF \ge f$ and $f(x_0)>0$, then $\norm{\wF}>0$ \rE{and further for any measurable function $u: \R_+ \to \R$  with $u(0)>0$}
\BQNY
 \lim_{\gamma \to \IF} P_{\gamma f - \widehat{\gamma f}}&=& \lim_{\gamma \to \IF} P_{\gamma f - \gamma \wwF}\\
&=& \lim_{\gamma \to \IF} P( B^H(t) + \gamma (f(t)- \wwF(t)) \le u(t), t\in \R_+) \\
&=& P( B^H(t)  \le u(t), t\in \R_+: f(t)= \wwF(t)) >0.
\EQNY
By Theorem \ref{ThA} for all  $\gamma$  large, 
\BQNY
P_{\gamma f} &\le & P_{\gamma f - \gamma \wF} \exp\Biggl(   - \frac{1}{2} \gamma^2 \norm{\widetilde{h}}^2 +
\gamma \int_0^\IF u(s) \, d (- K(s))\Biggr)\\
&=& P_{\gamma f - \gamma \wF} \exp\Biggl(   - \frac{1}{2} \gamma^2 \norm{\widetilde{h}}^2   (1+ o(1))\Biggr),
 \EQNY
hence  as $\gamma \to \IF$,  
\BQNY
 \ln P_{\gamma f} &\le &  - \frac{1}{2} \gamma^2 \norm{\widetilde{h}}^2  (1+ o(1)) +\ln P_{\gamma f - \gamma \wF}  
 =    - \frac{1}{2} \gamma^2 \norm{\widetilde{h}}^2 (1+ o(1)).                       
 \EQNY
It is clear that we can find $u_{-}$ such that $u_{-}(t)< u(t), t\in (0,\IF)$ such that $\int_{0}^\IF u_{-}(t) d (-K(t))$ is finite and
 $\pk{u_{-}(t) < B_H(t) \le u(t), t\in R_{+}}>0$. Applying again  Theorem \ref{ThA} for such $ u_{-}$ we have
\BQNY
\ln P_{\gamma f} & \ge  & -\frac{1}{2} \gamma^2 \norm{\widetilde{h}}^2 (1+o(1))
\EQNY
as $\gamma\to \IF$, and thus the claim follows.} \QED

{\bf Proof of Corollary \ref{crB}:}  In view of \eqref{LD} and the result of Corollary \ref{crA}, we have  
$$ \frac{1}{2} \gamma^2  \inf_{g\in \kH, g\ge f} \norm{g}^2_{ \Nhh } \sim \frac{1}{2} \gamma^2 \norm{\widetilde{h}}^2 $$
 as $\gamma \to \IF$. Since further $\norm{\wF}_{H}= \norm{\widetilde{h}}$ and the solution of the minimization problem is unique, then
$\wwF$ is its solution, thus the claim follows. \QED

\section{Appendix}
We present next one technical result.

\begin{lem} \label{AppE}  Let $H\in (1/2,1)$ {and suppose that the function} $g:\R_+\rightarrow \R_+$
{is such that} $g(t)=\int_0^t \rE{(K_{0+}^{H,*} f')}(s) ds$ 
for some $f$ such that $\rE{(K_{0+}^{H,*} f')}\in L_2(\R_+)$ and $f(0)=0$. Then $f(t)\geq 0, t\in \R_+$, holds.  

\end{lem}

{\bf Proof}:  We have that $$g(t)=\int_0^t\mathcal{D}_{0+}^{H-\frac12}(f'(u)u^{\frac12-H})(s)s^{H-\frac12}ds, \quad
\text{ with }f'(u)u^{\frac12-H}=I_{0+}^{H-\frac12}(g'(t)t^{\frac12-H})(u)$$
and
\begin{equation}\label{eq7.1}\begin{gathered}
f(u)=\int_0^us^{H-\frac12}I_{0+}^{H-\frac12}(g'(t)t^{\frac12-H})(s)ds\\
 =\Big(\Gamma \Big(H-\frac12\Big)\Big)^{-1}\int_0^u\Big(\int_s^uz^{H-\frac12}(z-s)^{H-\frac32}dz\Big) g'(s)s^{\frac12-H}ds.
\end{gathered}
\end{equation}

Setting  $r(s)=s^{\frac12-H}\int_s^uz^{H-\frac12}(z-s)^{H-\frac32}dz$, we may further write  
\BQNY
f(u)&=&-\Big(\Gamma \Big(H-\frac12\Big)\Big)^{-1}\int_0^ug(s)r'(s)ds
\EQNY
and
\BQNY
-r'(s)&=&-\Bigl(s^{\frac12-H}\int_s^uz^{H-\frac12}(z-s)^{H-\frac32}dz \Bigr)'_s =-(s^{\frac12-H}\int^{u-s}_0(z+s)^{H-\frac12}z^{H-\frac32}dz)'_s
\\
&=&(H-\frac12)s^{-\frac12-H}\int^{u-s}_0(z+s)^{H-\frac12}z^{H-\frac32}dz +s^{\frac12-H}u^{H-\frac12}(u-s)^{H-\frac32}
\\
&&-(H-\frac12)s^{\frac12-H} \int^{u-s}_0(z+s)^{H-\frac32}z^{H-\frac32}dz\\
&=&s^{\frac12-H}u^{H-\frac12}(u-s)^{H-\frac32}
+\cE{\int^{u-s}_0}(H-\frac12)s^{-\frac12-H}(z+s)^{H-\frac32}z^{H-\frac12}dz>0,
\EQNY
whence the  claim follows. \QED

\textbf{Acknowledgments.} E.\ Hashorva and Y. Mishura acknowledge support from
the Swiss National Science Foundation Grant 200021-1401633/1.  E.\ Hashorva kindly acknowledges partial support from the project RARE -318984, a Marie Curie IRSES Fellowship within the 7th European Community Framework.

\bibliographystyle{plain}
 \bibliography{boundaryCrossingsC}

\end{document}